\documentclass[leqno]{article}

\usepackage{amsmath}
\usepackage{amssymb}
\usepackage{amsthm}
\usepackage{enumitem}
\usepackage[T1]{fontenc}
\usepackage{hyperref} 

\title{Ueber die algebraischen Curven von den Geschlechtern $p = 4, \, 5$ und $6$, welche eindeutige Transformationen in sich besitze }
\author{Written by A. Wiman. \\
Translated by Linden Disney-Hogg, \\
Andrew Beckett, and Isabella Deutsch.}
\date{09/10/1895}

\begin{document}

\maketitle

\section*{Preface from the translators}
This is a translation of the Wiman paper,
\begin{center}
	Wiman, Anders. ``Ueber die algebraischen Curven von den Geschlechtern $p= 4, \, 5$ und $6$, welche eindeutigen Transformationen in sich besitzen."  Bihang till Kongl. Svenska vetenskaps-akademiens handlingar 21 (1895): 1-41.
\end{center}
which translates as 
\begin{center}
	On the algebraic curves of genus $p=4, \, 5$ and $6$, which possess unambiguous transformations into themselves, 
\end{center}
an influential paper written in German. \\
Care was taken to try and faithfully recreate the formatting of the original article, but this was not always possible. In contrast to the original, footnotes are numbered continuously throughout the document. Margin notes are included to anchor the pages of the original document to this translation. \\
LDH performed the original transcription and translation of the paper, and maintains the \href{https://github.com/DisneyHogg/Wiman-Translation}{GitHub copy} where errata may be reported. AB and ID contributed equally in refining the translation and providing proofreading. \\
The translators would like to thank the Royal Swedish Academy of Sciences and the Institut Mittag-Leffler for their permission to produce this English translation. This work is licensed under a Creative Commons Attribution-ShareAlike 4.0 International License.
\part*{The Paper}
\section*{Introduction}
The present treatise is intended to be a continuation of a work I have already published,\footnote{Bihang till K. Vet.Akad. Handlingar, Vol. 21, Section. I. No. 1.} which, among other things, treated the unambiguous correspondences on the non-hyperelliptic curves of genus $p = 3$. 
For the sake of convenience, we take the following equation from it, which was constructed using a theorem given by Mr. ZEUTHEN:
\begin{equation}\label{eq: A}
2(p-1) = 2n(p^\prime-1) + \sum \frac{n}{n_i}(n_i-1) \, . \tag{A}
\end{equation}
Here, on a curve of genus $p$, there should be an unambiguous correspondence of period $n$; the summation is over the coincidences where $n$ related points coincide in $\frac{n}{n_i}$ by $n_i$; finally, $p^\prime$ denotes the genus of a curve, the points of which are in a one-to-one relationship with the point groups of $n$ points that belong together through the correspondence.\footnote{Refer once again to a work by Mr. Hurwitz (Math. Ann., Vol. 32, p. 290).} \\
Since we refrain from the hyperelliptic case, which is easy to handle, we can bring the relevant structures to the projectively invariant normal curve $C_{2p-2}$ in the $ (p-1) $-dimensional space, whose coordinates are proportional to the ABELian integrands of the first kind. From this arises the advantage for us that the unambiguous correspondences are represented by collineations of the curve.\footnote{See the 5\textsuperscript{th} section of my previously-mentioned work.}
\section{p.4}\marginpar{start p.4}
For $ p = 4 $ we get a normal curve of the sixth order in ordinary three-dimensional space. This curve must be located on a surface of degree two $ F_2 $, which, however, can also degenerate into a cone $ K_2 $.\footnote{Noted in an essay by Mr. WEBER, Math. Ann., Vol. 13, p. 35.} The generators of this surface are formed by the trisecants of the curve. \\
CAYLEY has given a general theory for the singularities of a twisted curve,\footnote{See the related presentation by SALMON-FIEDLER, \textit{Analytische Geometrie des Raumes} II, 3. Edition, p. 105.}, and we do not find it inappropriate to quote some of his results for our case, since the characteristic numbers of the invariant normal curve also designate properties of the general algebraic structures of the relevant genus. \\
The $C_6$ mentioned has no true double points or cusps, since they can never occur on the $\varphi$-normal curve. The only possible relevant singularities are \textit{stationary tangents} $\Theta$ and double osculating planes $\Delta$; but those must be generators of the surface $ F_2 $, because every other degree has at most two points with $ F_2 $, and thus with the curve, in common. We shall find that these singularities appear in connection with the collineations of the curve. 
\begin{itemize}[label={}]
	\item Furthermore, let $ m $ be the  \textit{order} of the curve; 
	\item $ r $ the order of the developable surface formed by the tangents or the class of the perspective cone (the \textit{rank} of the curve);
	\item $ n $ the \textit{class} of the curve;
	\item $\alpha$ the number of \textit{stationary planes};
	\item $ h $ the number of \textit{apparent double points};
	\item $ g $ the number of lines in two osculating planes which lie in a given plane;
	\item $ x $ the order of the \textit{double curve};
	\item $ y $ the class of the double-touching developable surface.
\end{itemize}
Two systems of equations are obtained, the first by using PL\"UCKER's formulas for an arbitrary perspective cone over the given curve. For the order, class, double points, cusps, inflexions, double tangents one obtains the numbers $ m $, $ r $, $ h $, $ o $, $ n + \Theta $, $ y $ respectively. So given
\marginpar{start p.5}\[
m=6 \, , \quad p=4
\]
in advance, one finds,
\[
r = 18 \, , \quad h = 6 \, , \quad n = 36-\Theta \, , \quad y = 96 \, .
\]
So consider the cross-sectional curve of the developable surface. Here are the PL\"UCKER character numbers: $r$, $n$, $x$, $m+\Theta$, $\alpha$, $g+\Delta$. Thus 
\[
x = 126-\Theta \, , \quad \alpha=60-2\Theta \, , \quad g = 531-\frac{65}{2}\Theta + \frac{\Theta^2}{2} - \Delta \, . 
\]
The points of contact of the $ \alpha$-planes have a special meaning as those points of the structure of genus $p$ (4), in which a unique algebraic function of lower than $(p + 1)$th (5th) order becomes infinite without becoming infinite in other places. These places have been considered by many authors and given for their number $ (p-1) p (p + 1) $, so in our case 60; this is also confirmed if only the $ \Theta$-points are calculated twice. \\
The developable surface of order 18 and surface $ F_2 $ intersect in a curve of order 36; this can only be calculated twice from the cuspidal curve and consist of a number of common generators. Each generating system of $ F_2 $ delivers 12 of these, which can also be seen from the fact that each $ F_2$-generator has to divide the developable into 18 points, namely three double points on the curve and twelve common generators of the other system. Two of these 24 common straight lines move together in a stationary tangent that may appear, which gives the upper value 12 for $ \Theta $, which is actually reached. \\
\textit{Meanwhile we have found 2 groups of points with projective invariant properties, which under a collineation of the curve can only be interchanged with each other: the 60 contact points of the stationary planes and the 24 contact points of the common generators.} In addition, one could take the 24 points that cut the latter straight line from the curve. But this does not give a new group in two actually occurring cases: if all common generators are stationary tangents \marginpar{start p.6}or if all tangents in the above-mentioned cut-out points yield generators of the other $F_2$-system. \\
CAYLEY also treated other singularities by means of the theory of reciprocal surfaces.\footnote{See, for example, the presentation given by SALMON-FIELDER (Ibid p. 660).} Let us only give the following results. Let:
\begin{itemize}[label={}]
	\item $ \gamma^\prime $ be the number of osculating planes which also contain a tangent;
	\item $ t $ the number of points through which three tangents pass;
	\item $ t^\prime $ the number of planes containing three tangents.
\end{itemize}
The following formulas apply here if, of the relevant singularities, only $\Theta$ and $\Delta$ occur:
\begin{align*}
	\gamma^\prime &= rn + 12r - 14n - 6m - 8\Theta - 4\Delta  \, ; \\
	t &= \frac{1}{6}\left(r^3 - 3r^2 - 58r - 3r \left(n + 3m + 3\Theta \right) +42n + 78m + 78\Theta \right) \, ; \\
	t^\prime &= \frac{1}{6} \left(r^3 - 3r^2 - 58r - 3r \left(m + 3n + 3\Theta \right)+42m + 78n + 78\Theta\right) \, , 
\end{align*}
or in our case:
\[
\gamma^\prime = 324 - 12\Theta - 4\Delta \, ; \quad t=480-12\Theta \, ; \quad t^\prime = 120 \, . 
\]
Here we have found new invariant point groups. The $ t^\prime $ in particular have attracted attention in the general theory of curves, namely in the question of the adjoint $ \varphi$-curves, which each touch the base curve simply at $ p-1 $ points.\footnote{See CLEBSCH-LINDEMANN, \textit{Vorlesungen \"uber Geometrie} I, p. 847.} For the number of solutions one has $ 2^{p-1} (2^p-1) $, which for $ p = 4 $ agrees with the one already obtained.
\section{p.6}
The collineations of $C_6$ into itself must also transform the $F_2$ into itself. The two systems of generators can either be transformed into themselves or exchanged with one another. A special transformation can also occur which leaves all generators of a species individually fixed; two points on each are fixed, which together satisfy two fixed generators of the other kind, and the transformation has a period of three, because the points of the curve on the fixed generators must be cyclically exchanged; the six points on the two generators\marginpar{start p.7} of the other kind mentioned remain fixed and must have stationary tangents. \\
It can be seen from this that \textit{no new simple group can appear in the transformation groups of a structure of genus $ p = 4 $. Either the generating systems can be exchanged, and then half of the operations yield a normal subgroup which transforms each into itself, or, if exchange does not occur, the group is holohedrally isomorphic with one of the known finite binary groups, unless that every generator of a kind is transformed into itself by a cyclic $ G_3 $, but where this $ G_3 $ must have the normal property.} (The latter also applies to the case where the surface $ F_2 $ merges into a cone, and so the two systems coincide.)\\
If one uniquely assigns the ratios $ x_1: x_2 $ or $ y_1: y_2 $ to the two sets to be generated, then it is initially evident that $ C_6 $ can be found by a homogeneous equation $ F (x_1, x_2, y_1, y_2) = 0 $, which is degree 3 in both $ x_i $ and $ y_i $. This parameter setting is achieved by choosing $ xz = yw $ for the equation of $F_2$  and then using
\begin{equation}\label{eq: 2.1}
x: y : z : w = x_2y_1 :x_1 y_1 : x_1 y_2: x_2 y_2 \tag{1}
\end{equation}
to determine the ratios $ x_1: x_2 $ and $ y_1: y_2 $. From the point $ x = y = z = 0 $ or $ x_1 = y_1 = 0 $ one can now project the area and thus the curve onto the plane $ w = 0 $, whereby the image curve must satisfy the equation $F(y, x, y, z) = 0$.\footnote{See for example SALMON-FIELDER Ibid p. 521 or CLEBSCH-LINDEMANN, Vorlesungen \"uber Geometrie II 1, p. 422.} \\
A collineation of the first kind can be replaced with two binary substitutions as far as the area is concerned:
\begin{align*}
	x_1^\prime &= a_{11} x_1 + a_{12} x_2 \, , \quad x_2^\prime = a_{21} x_1 + a_{22} x_2 \, ; \\
	y_1^\prime &= b_{11} y_1 + b_{12} y_2 \, , \quad y_2^\prime = b_{21} y_1 + b_{22} y_2 \, . 
\end{align*}
Two generators of each kind merge into one another, and their 4 points of intersection remain fixed, but if one substitution is identical, all points of the fixed generators\marginpar{start p.8} of the other kind remain fixed. If one chooses the basic elements,\footnote{It is well known that coincident basic elements do not occur in finite groups.} one gets the simpler form:
\begin{equation}\label{eq: 2.2}
x_1^\prime = a_1 x_1 \, , \quad x_2^\prime = a_2 x_2 \, , \quad y_1^\prime = b_1 y_1 \, , \quad y_2^\prime = b_2 y_2 \, . \tag{2}
\end{equation}
The collineation in space is now determined from \eqref{eq: 2.1} and \eqref{eq: 2.2}:
\begin{equation}\label{eq: 2.3}
x^\prime = a_2 b_1 x_1 \, , \quad y^\prime = a_1 b_1 y \, , \quad z^\prime = a_1 b_2 z \, , \quad w^\prime = a_2 b_2 w \, . \tag{3}
\end{equation}
In a collineation of the second kind, where the generating sets are exchanged, one can initially start from substitutions of the following form:
\begin{align*}
	x_1^\prime &= \alpha_{11} y_1 + \alpha_{12} y_2 \, , \quad x_2^\prime = \alpha_{21} y_1 + \alpha_{22} y_2 \, ; \\
	y_1^\prime &= \beta_{11} x_1 + \beta_{12} x_2 \, , \quad y_2^\prime = \beta_{21} x_1 + \beta_{22} x_2 \, . 
\end{align*}
It is easy to see that through the requirements:
\[
x_1^\prime : x_2^\prime = x_1 : x_2 \, , \quad y_1^\prime: y_2^\prime = y_1 : y_2 \, , 
\]
determine two points fixed in the collineation on the surface, provided these conditions are not met in an infinite number of points. We choose generators through fixed points as basic elements and then get collineations: 
\begin{equation}\label{eq: 2.4}
x_1^\prime = \alpha_1 y_1 \, , \quad x_2^\prime = \alpha_2 y_2 \, , \quad y_1^\prime = \beta_1 x_1 \, , \quad y_2^\prime = \beta_2 x_2 \, . \tag{4}
\end{equation}
The collineation of space now results from \eqref{eq: 2.1} and \eqref{eq: 2.4}: 
\begin{equation}\label{eq: 2.5}
x^\prime = \alpha_2 \beta_1 z \, , \quad y^\prime = \alpha_1 \beta_1 y \, , \quad z^\prime = \alpha_1 \beta_2 x \, , \quad w^\prime = \alpha_2 \beta_2 w \, . \tag{5}
\end{equation}
If this operation is repeated, a collineation of the first kind is obtained:
\begin{equation}\label{eq: 2.6}
x^\prime = \alpha_1 \alpha_2 \beta_1 \beta_2 x \, , \quad y^\prime = \alpha_1^2 \beta_1^2 y \, , \quad z^\prime = \alpha_1 \beta_1 \alpha_2 \beta_2 z \, , \quad w^\prime = \alpha_2^2 \beta_2^2 w \, , \tag{6}
\end{equation}
but not an arbitrary one, since all points of the straight line $ y = w = 0 $ remain fixed. \\
If we also emphasise that the equation \eqref{eq: A} for $ p = 4 $ only allows the prime number solutions 2, 3, 5, then we now have the means to identify the possible collineation groups of a structure of genus $ p = 4 $ as quickly as possible.
\section{p.9}\marginpar{start p.9}
We now want to visit these collineations. The equations of the structures in question are written with the smallest possible number of constants, which means that the arbitrary moduli are also given. \\
For a cyclic group of period 2, either the generators can be swapped or not. In the first case one gets an infinite number of fixed points, the intersections of two mutually exchanging generators. These form a conic section, where the pole of the plane of the conic section with respect to the $F_2$ provides the centre of perspective. The equation of the curve is obtained by choosing two fixed points on the curve for $x_1 = y_1 = 0$ and $x_2 = y_2 = 0$ and the simultaneous interchangeability of $x_1$, $y_1$ and $x_2 $, $y_2$ takes into account the following:
\begin{align}\label{eq: 3.1}
 & x_1^2 y_1^2 (x_1 y_2 + x_2 y_1) + x_1 y_1(x_1^2 y_2^2 + x_2^2 y_1^2 + a x_1 x_2 y_1 y_2) + \nonumber \\
 +& b(x_1^3 y_2^3 + x_2^3 y_1^3) + c x_1 x_2 y_1 y_2 (x_1 y_2 + x_2 y_1) +  \tag{1}\\
 +& x_2 y_2 (d(x_1^2 y_2^2 + x_2^2 y_1^2) + ex_1 x_2 y_1 y_2) + f x_2^2 y_2^2  (x_1 y_2 + x_2 y_1) =0 \, . \nonumber
\end{align}
If we go to the coordinates of space, we find (according to Eq. \eqref{eq: 2.5} \S 2, where $\alpha_1 = \alpha_2 = \beta_1 = \beta_2$) $x-z=0$ for the perspective plane  and $y=w=x+z=0$ for the perspective centre. The curve must be on a cone of the 3rd order, whose vertex is the centre, so that every two points on the same generator are interchanged by the collineation. The tangents of the 6 points in the perspective plane pass through the centre and the associated planes are stationary. The planes of the 9 turning points of the cone are doubly osculating, thus forming 9 $\Delta$-planes. \\
Now the cyclic $G_2$ may not effect any permutation of the generator sets. One can bring it to the form:
\[
x_1^\prime = -x_1 \, , \quad y_1^\prime = - y_1 \, , \quad x_2^\prime = x_2 \, , \quad y_2^\prime = y_2 \, .   
\]
The equation can be written in two equivalent ways, from which we choose the following:
\begin{align}\label{eq: 3.2}
	x_1^3 y_1^3 + x_1 y_1 x_2 y_2 (ax_1 y_1 + b x_2 y_2) +& c x_2^3 y_2^3 +x_1^3 y_1 y_2^2 + x_1 x_2^2 y_1^3 + \tag{2} \\ 
	+& d x_2^3 y_1^2 y_2 + e x_1^2 x_2 y_2^3 = 0 \, .  \nonumber 
\end{align}
\marginpar{start p.10}Two points of the curve are transformed into themselves: $x_1=y_2=0$ and $x_2=y_1=0$. If one takes account of the equations \eqref{eq: 2.2} and \eqref{eq: 2.3} of the previous section, one finds that in this collineation every point of the straight lines $x=z=0$ and $y=w=0$ remains fixed. Any straight line that intersects these two lines must merge into itself. The connecting lines of corresponding points on the curve thus form a ruled surface, for which those two straight lines are guidelines, and it is easily seen that this ruled surface is of fifth order and contains $y=w=0$ as triple $x=z=0 $ as a double directrix. \\
The possible location of the $C_6$ on a 3rd order cone or a 5th order ruled surface thus characterises the two different types of involutory-unambiguous correspondences on the curve. Since only 6 or 5 moduli are included in \eqref{eq: 3.1} and \eqref{eq: 3.2} and 9 moduli are available for the general $C_6$, 3 or 4 conditions are fulfilled there. The cone $K_3$, or the ruled surface $R_5$, designates a structure of the genus $p^\prime=1$, or 2, to which the $C_6$ is 1-2-unambiguously related. No other rectilinear surfaces can simply contain the curve if at the same time the generators are supposed to be bisecants. \\
In the case of equation \eqref{eq: 3.2} we examine the bisecant ruled surface whose directrix is $y=w=0$, which itself contains 2 points of the curve. From each point of this straight line there are 5 bisecants (except $y=w=0$ itself) and in each plane through it there are 6, so that the order 11 is obtained. The ruled surface thus consists of the mentioned $R_5$ and a $R_6$, for which the $C_6$ is a double curve. The generators of this $R_6$, which start from a point of the directrix, merge into each other through the transformation $G_2$; but since every plane through the straight line remains invariant, those generators must also lie in the same plane through this straight line. One can therefore call $y=w=0$ a double tangent line for the $R_6$, since the cross-sectional curves there receive tangent nodes.\footnote{I have dealt with this $R_6$ in my writing \guillemotright Klassifikation af regelytorna af 6. graden\guillemotright{} p. 58 (Diss. Lund 1892)} This $R_6$ is usually of genus $p^\prime=2$, but splits into 2 $K_3$ when 
\begin{equation}\label{eq: 3.3}
d=e \, . \tag{3}
\end{equation}
\marginpar{start p.11}This equation says that the two points of $C_6$ on $y=w=0$ have the same plane $y+dw=0$. The transformation group of the curve into itself is in this case a \textit{Klein four-group}. Two $G_2$ interchanging the generating sets appear:
\[
x_1^\prime = \pm y_1 \, , \quad y_1^\prime = \pm x_1 \, , \quad x_2^\prime = y_2 \, , \quad y_2^\prime = x_2 \, .
\]
The associated perspective planes are $x \mp z = 0$ and the perspective centres $y=w = x \pm z=0$. \\
However, one can also obtain a \textit{Klein four-group of a different kind}, which leaves the generating groups uninterchanged altogether. We start again from Equation \eqref{eq: 3.2}, but at first we think of the ratios of all 8 coefficients as completely arbitrary. Under the remaining operations of the Klein four-group, the two generators of a family that are fixed under $G_2$ must be exchanged with each other. By means of suitable normalization, the two new $G_2$ can be given the following form:
\[
x_1^\prime = \pm x_2 \, , \quad x_2^\prime = x_1 \, , \quad y_1^\prime = \pm y_2 \, , \quad y_2^\prime = y_1 \, .
\]
The equation for the curve is:
\begin{align}\label{eq: 3.4}
	x_1^3 y_1^3 +& x_2^3 y_2^3 + ax_1 y_1 x_2 y_2 (x_1 y_1 + x_2 y_2) + b y_1 y_2 (x_1^3 y_2 + x_2^3 y_1) + \tag{4} \\ 
	+& c x_1 x_2 (x_1 y_2^3 + x_2 y_1^3) = 0 \, .  \nonumber 
\end{align}
A cyclic $G_4$ must have a $G_2$ as a subgroup, which leaves the generating groups uninterchanged. A look back at equation \eqref{eq: A} teaches us that the two fixed points of $C_6$ under $G_2$ must also maintain their fixed positions under $G_4$. In the collineations of period 4, the systems must be exchanged; this is already evident from the fact that the necessary number of points for the purpose of cyclic permutation are missing on the fixed generators occurring in this case. So we start from equation \eqref{eq: 3.2} under the most general assumption possible with regard to the coefficients. A collineation of period 4, which swaps the systems and keeps the points $x_1 = y_2 =0$ and $x_2 = y_1 = 0$ fixed, can be brought to the following form:
\[
x_1^\prime = y_2 \, , \quad y_2^\prime = -ix_1 \, , \quad x_2^\prime = iy_1 \, , \quad y_1^\prime = x_2 \, .
\]
\marginpar{start p.12} The generator for which $x_1=x_2$ can be chosen arbitrarily. We can therefore state that in the equation we are looking for, for example the coefficients for $x_1^3 y_1 y_2^2$ and $x_2^3 y_1^2 y_2$ become equal. Then we get equation \eqref{eq: 3.4} again, only with the specialisation that
\begin{equation}\label{eq: 3.5}
c = -b \, . \tag{5} 
\end{equation}
However, the curve determined in this way also contains the two groups of four belonging to cases \eqref{eq: 3.3} and \eqref{eq: 3.4}, so that the complete collineation group is of the nature of a \textit{dihedral} $G_8$. \\
One now easily proves that the 5 groups discussed are the only possible ones without collineations of higher prime number period.
\section{p.12}
We now add the cyclic groups of period 3. With such a group, either each generator of a group can merge into itself or not. In the latter case, the collineation can be created in the following form: 
\[
x_1^\prime = j^2 x_1 \, , \quad x_2^\prime = x_2 \, , \quad y_1^\prime = j y_1 \, , \quad y_2^\prime = y_2 \, , \quad \left( j = e^{\frac{2 \pi i}{3}} \right)
\] 
or in the coordinates of space (according to \eqref{eq: 2.2} and \eqref{eq: 2.3} \S 2):
\[
x^\prime = jx \, , \quad y^\prime = y \, , \quad z^\prime = j^2 z \, , \quad w^\prime = w \, . 
\]
The equation of a $C_6$ that merges into itself through this transformation can be written with only 3 moduli:
\begin{equation}\label{eq: 4.6}
x_1^3 y_1^3 + x_2^3 y_2^3 + a(x_1^3 y_2^3 + x_2^3 y_1^3) + b x_1^2 y_1^2 x_2 y_2 + c x_1 y_1 x_2^2 y_2^2 = 0 \, , \tag{6}
\end{equation}
whereby, of course, one disregards structures with double points and thus lower genus. But the curve \eqref{eq: 4.6} also has three perspective $G_2$ into itself: 
\[      
x_1^\prime = y_1 (j y_1 , j^2 y_1) \, , \quad y_1^\prime = x_1 \, , \quad x_2^\prime = y_2 \, , \quad y_2^\prime = x_2 (j x_2, j^2 x_2) \, .  
\]
or (according to \eqref{eq: 2.4} and \eqref{eq: 2.5} \S 2):
\[
x^\prime = z(j^2z, jz) \, , \quad y^\prime = y \, , \quad z^\prime = x(jx, j^2x) \, , \quad w^\prime = w \, . 
\]
The relevant perspective planes are: 
\[
x-z=0 \, , \quad x-j^2z=0 \, , \quad x-jz=0 \, , 
\] 
\marginpar{start p.13}and the perspective centres are obtained as the intersections of the straight lines $y=w=0$ with $x+z=0$, $x+j^2z=0$ and $x+jz=0 $. \\
The transformation group of the curve \eqref{eq: 4.6} is thus a \textit{dihedral group} $G_6$. Six related points lie in a plane through $y=w=0$. The straight lines connecting the corresponding points at the $G_2$ form 3 $K_3$, the vertices of which are the above-mentioned perspective centres. These $K_3$ are cyclically exchanged by the $G_3$. One can now easily deduce that a plane through $y=w=0$ must either touch all or none of $K_3$. There are 6 planes of that kind, which must also touch the curve three times; the points of contact belong to stationary planes, and their tangents each pass through one of the three vertices. \\
We now think of the $G_3$ as a subgroup of a higher group that does not swap systems. This group cannot be a cyclic $G_6$, since the non-fixed points permute each other cyclically to six on such a fixed generator, and not an octahedral group, because the cyclic $G_4$ required for this do not occur according to the previous section. The possible cases to be discussed here are thus dihedral $G_6$ and the tetrahedral group. \\
In the case of the other collineations of a linear dihedral group, those elements which remain fixed in the cyclic subgroup must be exchanged. In order for this exchange to be possible, a condition must be fulfilled in equation \eqref{eq: 4.6}, for which we can set 
\begin{equation}\label{eq: 4.7}
b=c \, . \tag{7}
\end{equation}
The complete collineation group of the curve is then of the type of a \textit{dihedral} $G_{12}$. As subgroups we mention the dihedral $G_6$, which has the general curve \eqref{eq: 4.6}, the dihedral $G_6$, which does not commute the generators, whose three collineations of period 2 are the following: 
\[
x_1^\prime = x_2( jx_2, j^2 x_2) \, , \quad x_2^\prime = x_1 \, , \quad y_1^\prime = y_2(j^2 y_2, j y_2) \, , \quad y_2^\prime = y_1 \, , 
\] 
and finally the cyclic $G_6$, which is formed by repeating the operation: 
\[
x_1^\prime = j y_2 \, , \quad y_2^\prime = x_1 \, , \quad x_2^\prime = y_1 \, , \quad y_1^\prime = j^2 x_2 \, . 
\]
\marginpar{start p.14}We now assume that the subgroup which does not interchange the systems is a tetrahedral group. As a normal subgroup, we have a Klein four-group here. We are therefore looking for the further conditions which the coefficients of equation \eqref{eq: 3.4} must satisfy in this case. The fixed elements in the individual operations of the Klein four-group are to be cyclically exchanged in the collineations of period 3 of the tetrahedron group. This collineation of period 3 can now easily be represented in the following way:
\[
x_1^\prime = \pm i(x_1 \pm x_2) \, , \quad x_2^\prime = x_1 \pm x_2 \, , \quad y_1^\prime = \pm i (y_1 \pm y_2) \, , \quad y_2^\prime = y_1 \pm y_2 \, , 
\]
and thus obtains for the equation of the curve merging into itself:
\begin{equation}\label{eq: 4.8}
(x_1 y_1 + x_2 y_2)^3 + b(x_1 y_1 - x_2 y_2) (x_1 y_2 - x_2 y_1)( x_1 y_2 + x_2 y_2) = 0 \, . \tag{8}
\end{equation}
The complete collineation group of curve \eqref{eq: 4.8} in itself turns out to be holohedrally isomorphic with an \textit{octahedral group}. The plane $x_1 y_1 + x_2 y_2 = 0$ or $y+w=0$ always merges into itself. The pole of this plane in relation to the area $F_2$ must therefore remain fixed for all 24 collineations. The curve lies on 6 cones $K_3$ whose vertices are in the invariant plane.
\section{p.14}
We now assume that each generator of the $y$ family merges into itself under a $G_3$. The $G_3$ can be generated by the collineation:
\[
x_1^\prime = j x_1 \, , \quad x_2^\prime = x_2 \, , \quad y_1^\prime = y_1 \, , \quad y_2^\prime = y_2 \, ,
\]
and the equation of an associated curve is:
\begin{equation}\label{eq: 5.9}
x_1^3 f_3 (y_1, y_2) + x_2^3 \varphi_3 (y_1, y_2) = 0 \, . \tag{9}
\end{equation}
This curve has 3 independent moduli. \\
We now consider the case that $f_3$ and $\varphi_3$ contain a common transformation group within them. The same must of course also hold for their HESSian covariants, and we first assume that these do not coincide. There are two possibilities here: either the points of one covariant can be base points of a $G_2$, \marginpar{start p.15}which swaps those of the other, or the points of both covariants are swapped by a $G_2$, whose base points must be supplied by its simultaneous quadratic covariant; but that case turns out to be impossible, because in order for the three points $f_3=0$ (resp. $\varphi_3=0$) to merge into each other through a $G_2$, one of them must remain fixed and thus provide a base point. \\
But one also comes to higher groups if $f_3=0$ and $\varphi_3=0$ can be exchanged with each other. Because the 6 points $f_3 \varphi_3 = 0$ have to be ordered cyclically into one or more series, the period of a collineation can only be 6 or 2.
By the same way of reasoning for the HESSian covariants, one obtains 2 as the only possible period. The equation of the curve can be found in the following form:
\begin{align}\label{eq: 5.10}
 & x_1^3(y_1 - y_2)(y_1 - ay_2)(y_1 - b y_2) + \tag{10}\\ 
 +& x_2^3(y_1 + y_2)(y_1 + a y_2) (y_1 + b y_2) = 0 \, . \nonumber 
\end{align}
The transformation group is \textit{a dihedral} $G_6$ whose collineations of period 2 are the following:
\[
x_1^\prime = x_2 (j x_2, j^2 x_2) \, , \quad x_2^\prime = x_1 (j^2 x_1 , j x_1) \, , \quad y_1^\prime = y_1 \, , \quad y_2^\prime = - y_2 \, . 
\]
Since the points of the HESSian covariants can be paired in two ways, one has the possibility to consider that $f_3 = 0$ and $\varphi_3 = 0$ are also swapped in two ways. But then $f_3=0$ and $\varphi_3=0$ must each contain a transformation different from the identity, which can be obtained by combining two permutations. The equation of the curve can be brought into the following form:
\begin{equation}\label{eq: 5.11}
x_1^3 y_1 (y_1^2 + a y_2^2) + x_2^3 y_2 (a y_1^2 + y_2^2) = 0 \, . \tag{11}
\end{equation}
The collineation group of this curve is a \textit{dihedral} $G_{12}$. A generating operation of the cyclic subgroup $G_6$ is found in:
\[
x_1^\prime = -j x_1 \, , \quad x_2^\prime = x_2 \, , \quad y_1^\prime =-y_1 \, , \quad y_2^\prime = y_2 \, .
\]
The 6 remaining collineations are the following:
\[
x_1^\prime = \pm x_2(\pm j x_2 , \pm j^2 x_2) \, , \quad x_2^\prime = x_1 \, , \quad y_1^\prime = \pm y_2 \, , \quad y_2^\prime = y_1 \, . 
\]
\marginpar{start p.16}We now examine the case where the HESSian covariants of $f_3$ and $\varphi_3$ coincide. The equation of the curve can easily be brought into the form:
\begin{equation}\label{eq: 5.12}
x_1^3 y_1^3 + x_1^3 y_2^3 + x_2^3 y_1^3 + a^3 x_2^3 y_2^3 = 0 \, . \tag{12}
\end{equation}
This curve has 12 stationary tangents whose tangents are cut out in threes each by the generators $x_1=0$, $x_2=0$, $y_1=0$ and $y_2=0$. Here one has a $G_3$ for each system, in which each generator remains fixed:
\begin{align*}
	x_1^\prime &= j x_1 \, , \quad  x_2^\prime = x_2 \, , \quad y_1^\prime = \phantom{j} y_1 \, , \quad y_2^\prime = y_2 \, ; \\
	x_1^\prime &= \phantom{j} x_1 \, , \quad x_2^\prime = x_2 \, , \quad y_1^\prime = j y_1 \, , \quad y_2^\prime = y_2 \, .
\end{align*}
By combining these $G_3$ one gets a $G_9$, which consists of 4 $G_3$ and is normal in the overall group. If one adds a collineation that exchanges the generators $x_1 = 0$ and $x_2 = 0$ or $y_1 = 0$ and $y_2 = 0$: 
\[
x_1^\prime = a x_2 \, , \quad x_2^\prime = x_1 \, , \quad y_1^\prime = y_2 \, , \quad y_2^\prime = \frac{1}{a} y_1 \, , 
\]
and combines this with $G_6$, one gets 9 collineations from the period 2. So one has a $G_{18}$, which converts every generator system into itself and is obviously holohedrally isomorphic with the collineation group of a general planar $C_3$. There are also 18 collineations that swap systems. One gets these by combining the $G_{18}$ with any of them, like:
\[
x_1^\prime = y_1 \, , \quad y_1^\prime = x_1 \, , \quad x_2^\prime = y_2 \, , \quad y_2^\prime = x_2 \, . 
\]
Of these 18 new collineations, 12 have period 6 and 6 have period 2. \\
With the mentioned 36 collineations of the curve \eqref{eq: 5.12}, 2 straight lines remain fixed: the straight lines connecting the points $x_1 = y_1 = 0$ and $x_2 = y_2 = 0$ or $x_1 = y_2 = 0$ and $x_2 = y_1 = 0$, or (Eq. \eqref{eq: 2.1} \S 2) $x=z=0$ and $y=w=0$. The transformation groups of these straight lines are dihedral $G_6$ and thus obtained with $G_{36}$ from the combination of these two $G_6$, and by combining subgroups of $G_6$ subgroups of $G_{36}$ arise. The vertices of three $K_3$ containing the curve lie on each of the two fixed straight lines.\marginpar{start p.17} The mentioned dihedral $G_6$ arise by summarising all possible permutations of these vertices. The $G_{36}$ is therefore to be regarded as a subgroup of the group of all possible permutations of 6 objects, namely the above-mentioned $G_{18}$ that does not permute the systems is obtained as the subgroup of $G_{36}$, in which permutations occur in an even number. \\
The question arises whether the curve \eqref{eq: 5.12} can contain further collineations. The generator $x_1=0$ can only be exchanged with the other three generators containing $\Theta$-points: $x_2=0$, $y_1=0$, $y_2=0$. So if $x_1=0$ merges into itself under $n$ collineations, then the order of the group is $4n$. The number $n$ is obtained by combining $G_3$, which leaves every point of the imagined generator unchanged, and the common transformation group of $f_3$ and $\varphi_3$ in itself. However, the latter group is usually a $G_3$ if the HESSian covariants coincide, but a dihedral group of order 6 if $f_3$ and $\varphi_3$ are mutual third-degree covariants. In this case, where $a^3=-1$ and thus the equation of the curve: 
\begin{equation}\label{eq: 5.13}
x_1^3 y_1^3 + x_1^3 y_2^3 + x_2^3 y_1^3 - x_2^3 y_2^3 = 0 \, , \tag{13}
\end{equation} 36 new collineations appear. One gets this by combining the $G_{36}$ with any new one, e.g
\[
x_1^\prime = x_2 \, , \quad x_2^\prime = x_1 \, , \quad y_1^\prime = y_1 \, , \quad y_2^\prime = - y_2 \, . 
\]
Of the 36 collineations thus added, those which interchange the systems are of period 4; of the rest, 12 have period 6 and 6 have period 2. Those lines fixed under $G_{36}$ are swapped under the 36 new collineations, and one can consider the $G_{72}$ as a group of permutations of 6 objects, where about the three first and three last things always form the same 2 commuting groups. The $G_{72}$ has three normal subgroups of order 36: the collineation group of the general curve \eqref{eq: 5.12}, the subgroup that does not interchange the systems, and finally a group of even permutations of the 6 $K_3$ vertices. The latter is holohedrally isomorphic with the collineation group of the plane harmonic $C_3$. As mentioned, \marginpar{start p.18}the curve \eqref{eq: 5.12} lies on 6 $K_3$, the vertices of which are occupied three times each on the two distinguished straight lines. The same curve also lies on 9 $R_5$, whose triple directrix connects two cone vertices. The special curve \eqref{eq: 5.13} also lies on 6 other $R_5$ whose triple directrices each connect two $\Theta$-points; these lines are double tangent straight lines of 6 $R_6$ for which \eqref{eq: 5.13} is double curve. \\
In general, 72 points of the curve \eqref{eq: 5.13} permute in the $G_{72}$ as a closed system. But the 12 $\Theta$-points only permute among themselves, and the same applies to the system of the 36 $\alpha$-points, and finally also to 18 points in the 9 straight lines connecting two $K_3$ vertices, which form the points of contact of 9 double-osculating planes. Each of these points remains fixed under a $G_6$, or $G_2$ or $G_4$. 
\section{p.18}
It remains to consider the case where the curve is transformed into itself by a $G_5$. With an odd number of periods, of course, the generating systems cannot be interchanged. Since the points of the curve on the 4 fixed generators cannot be cyclically permuted to 5 each, they must be occupied in the 4 fixed points of intersection. However, this can only be achieved by simply touching the curve in each of the 4 points. It is now easy to see through experiments that the equation of a curve of the type considered here with a collineation of period 5 can be brought into the following form: 
\begin{equation}\label{eq: 6.14}
x_1^3 y_1^2 y_2 + x_1^2 x_2 y_2^3 + x_1 x_2 y_2^3 + x_1 x_2^2 y_1^3 + a^5 x_2^3 y_1 y_2^2 = 0 \, , \tag{14}
\end{equation}
where if $\varepsilon^5=1 $, the collineations are represented as follows:
\[
x_1^\prime = \varepsilon x_1 \, , \quad x_2^\prime = \varepsilon^4 x_2 \, , \quad y_1^\prime = \varepsilon^2 y_1 \, , \quad y_2^\prime = \varepsilon^3 y_2 \, . 
\]
However, the curve \eqref{eq: 6.14} also has 5 collineations of period 2: 
\[
x_1^\prime = a^2 \varepsilon^4 x_2 \, , \quad x_2^\prime = \frac{\varepsilon}{a} x_1 \, , y_1^\prime = \varepsilon^3 y_2 \, , \quad y_2^\prime = \frac{\varepsilon^2}{a} y_1 \, , 
\]
so that the complete group of its collineations is a \textit{dihedral} $G_{10}$. \marginpar{start p.19}One easily finds that no collineation of period 10 can transform the curve \eqref{eq: 6.14} into itself; because either all 4 fixed points of $G_5$ would also remain fixed here, which contradicts equation \eqref{eq: A}, or 2 of them would be swapped with each other, which would be geometrically inconsistent. Collineations with a higher number of periods are also impossible for the same reasons. One only has to consider the case in which a generator is transformed into itself by an icosahedron group. Because the $G_3$ occurring here must be of the type dealt with in the 4th section, and the systems were always exchanged there, the complete group must be of order 120. With such a curve with 120 transformations in itself, the 24 tangents, which are also generators of $F_2$, form 6 quadrilaterals corresponding to the 6 $G_5$ within the icosahedron group, whose corners lie in the points of contact. \\
If an exchange of the systems for the curve \eqref{eq: 6.14} is to be possible, the quadrilateral on which it is based must be converted into itself if one does not want to assume several such quadrilaterals, and thus also several $G_5$ and one $G_{120}$. One easily finds the condition $a^5=-1$; then there are 5 cyclic $G_4$ whose generating operations are the following: 
\[
x_1^\prime = \varepsilon^2 y_1 \, , \quad y_1^\prime = - \varepsilon^4 x_2 \, , \quad x_2^\prime = \varepsilon^3 y_2 \, , \quad y_2^\prime = \varepsilon x_1 \, , 
\] 
where $\varepsilon$ is a fifth root of unity. \\
Meanwhile Mr. GORDAN recognised\footnote{See P. GORDAN, \textit{\"uber die Aufl\"osung der Gleichungen vom f\"unften Grade} (Math. Ann. Vol. 13) or F. KLEIN, \textit{Vorlesungen \"uber das Ikosaeder}, p. 194ff.} that the curve obtained in this way: 
\begin{equation}\label{eq: 5.15}
x_1^3 y_1^2 y_2 + x_1^2 x_2 y_2^3 + x_1 x_2^2 y_1^3 - x_2^3 y_1 y_2^2 = 0 \tag{15}
\end{equation} 
allows 120 collineations. Mr. KLEIN calls this curve \textit{Bring's curve}, for the following reason. If one introduces pentahedron coordinates $z_1$, $z_2$, $z_3$, $z_4$, $z_5$ with the identical relation 
\[
\sum z_i = 0 \, ,
\]
\marginpar{start p.20}one obtains the same as an intersection curve of two surfaces, the main surface and the diagonal surface: 
\begin{equation}\label{eq: 5.15a}
\sum z_i^2 = 0 \, , \quad \sum z_i^3 = 0\footnote{One can consider the coordinates of a point on the curve as the roots of an equation of the 5th degree, which has been brought into BRING's normal form by Tschirnhaus transformation.} \, . \tag{15a}
\end{equation}
The 120 collineations of the curve result from the permutations of the $z_i$, and the 60 even permutations leave the generating systems unchanged. \\
The nature of the individual collineations is briefly indicated here. By cyclic permutation of all $z_i$ one gets 24 collineations of period 5. Twenty collineations of the type:
\[
z_1^\prime = z_2 \, , \quad z_2^\prime = z_3 \, , \quad z_3^\prime = z_1 \, , \quad z_4^\prime = z_5 \, , \quad z_5^\prime = z_4 \, ,
\] 
have the period 6. The resulting 10 cyclic $G_6$ have subgroups of period 3, which yield 20 new collineations, E.g.: 
\[
z_1^\prime = z_2 \, , \quad z_2^\prime = z_3 \, , \quad z_3^\prime = z_1 \, , \quad z_4^\prime = z_4 \, , \quad z_5^\prime = z_5 \, ;
\] 
In addition, there are the subgroups of period 2:
\[
z_1^\prime = z_2 \, , \quad z_2^\prime = z_2 \, , \quad z_3^\prime = z_3 \, , \quad z_4^\prime = z_5 \, , \quad z_5^\prime = z_4 \, . 
\]
These 10 collineations are perspective; the perspective plane in the chosen example is $z_4 - z_3 = 0$ and the perspective centre $z_1 = z_2 = z_3 = z_4+z_5$. The 60 points of contact of the stationary planes are in the perspective planes. The 10 straight lines $z_i=z_k=0$ each go through three perspective centres; then the straight line $z_1 = z_2 = 0$ goes through the centres assigned to the planes $z_3=z_4$, $z_4=z_5$ and $z_5 = z_3$ and stays with the permutations of $z_3$, $z_4$, $z_5$ resulting dihedral $G_6$. Thirty collineations have period 4; in these, one pentahedron plane remains fixed and the 4 others are cyclically interchanged. By repeating one of these, 15 collineations of period 2 \guillemotright with axes\guillemotright{} arise. As an example we choose:
\[
z_1^\prime = z_2 \, , \quad z_2^\prime = z_1 \, , \quad z_3^\prime = z_4 \, , \quad z_4^\prime = z_3 \, , \quad z_5^\prime = z_5 \, . 
\]
\marginpar{start p.21}The one axis $z_1 + z_2 = z_3 + z_4 = z_5 = 0$ goes through 2 of the aforementioned perspective centres: $z_1 +z_2 = z_3 = z_4 = z_5 = 0$ and $z_3 + z_4 = z_1 = z_2 = z_5=0$, and in fact each of the 15 connecting lines is such an axis. According to the discussion in Equation \eqref{eq: 3.2} of section 3, each such axis passes through the tangents of a doubly osculating plane. \\
The higher subgroups within the $G_{120}$ are: the icosahedron groups, which leave the systems uninterchanged, 5 octahedron groups, in each of which one pentahedron plane remains fixed, 6 $G_{20}$ of the kind previously discussed, 10 dihedral $G_{ 12}$, in which the pentahedron planes permute in 2 systems of 2 and 3, and finally dihedral groups of the order 10, 8, 6 and 4, which, however, form only subgroups of the subgroups discussed. \\
Corresponding to the 10 perspective $G_2$ and the 15 $G_2$ \guillemotright with axes\guillemotright, BRING's curve lies at 10 $K_3$ and 15 $R_5$. Each turning plane of this $K_3$ osculates the curve twice. However, 15 turning points, the mentioned axes, belong to 2 $K_3$. Sixty others are obtained, so that the total number of planes doubly osculating the curve is 75. Of those planes $\gamma^\prime$ which touch the curve at one point and osculate at another, 300 are absorbed (according to a formula at the end of the 1st section) by the 75 $\Delta$-planes; the remaining 24 osculate in the points whose tangents are generators of $F_2$, because the other generator through such a point also contains a tangent of the curve. One could easily draw out several such remarkable properties of BRING's curve. \\
The 60 $\alpha$-points, the 30 points of tangency of the system of 15 $\Delta$-planes and the 24 points in which the tangents are generators of the $F_2$ form the only systems of less than 120 points that are permuted in a closed manner under $G_{120}$ .
\section{p.21}
Now that we have found all possible collineation groups of a structure of genus $p=4$ whose normal curve of the 6th order lies on a genuine surface of the 2nd degree, the results obtained are summarised here. $M$ denotes the number of independent moduli of the entities concerned. \marginpar{start p.22}
\begin{enumerate}[label=\alph*.]
	\item $M=6$.
	\begin{enumerate}[label=\arabic*)]
		\item $G_2$, which swaps the generating systems of the surface $F_2$. Equation \eqref{eq: 3.1}.
	\end{enumerate}
    \item $M=5$.
    \begin{enumerate}[label=\arabic*)]
    	\setcounter{enumii}{1}
    	\item $G_2$ without swapping the systems. Eq. \eqref{eq: 3.2}.
    \end{enumerate}
    \item $M=4$.
    \begin{enumerate}[label=\arabic*)]
    	\setcounter{enumii}{2}
    	\item Klein four-group with swapping of the systems. Eq. \eqref{eq: 3.3}.
    \end{enumerate}
    \item $M=3$.
    \begin{enumerate}[label=\arabic*)]
    	\setcounter{enumii}{3}
    	\item Klein four-group without swapping of the systems. Eq. \eqref{eq: 3.4}.
    	\item Dihedral $G_6$ with exchange of systems. Eg. \eqref{eq: 4.6}. 
    	\item $G_3$. Eq. \eqref{eq: 5.9}.
    \end{enumerate}
    \item $M=2$.
    \begin{enumerate}[label=\arabic*)]
    	\setcounter{enumii}{6}
    	\item Dihedral $G_8$ with exchange of systems. Eq. \eqref{eq: 3.5}.
    	\item Dihedral $G_{12}$ with exchange of systems. Eq. \eqref{eq: 4.7}.  
    	\item Dihedral $G_6$ without swapping of the systems. Eq. \eqref{eq: 5.10}.
    \end{enumerate}
    \item $M=1$. 
    \begin{enumerate}[label=\arabic*)]
    	\setcounter{enumii}{9}
    	\item Octahedral group. Eq. \eqref{eq: 4.8}. 
    	\item Dihedral $G_6$ without swapping of the systems. Eq. \eqref{eq: 5.11}.
    	\item Dihedral group of order 10 without interchanging the systems. Eq. \eqref{eq: 6.14}.
    	\item Group of Order 36. Eq. \eqref{eq: 5.12}.
    \end{enumerate}
    \item $M=0$. 
    \begin{enumerate}[label=\arabic*)]
    	\setcounter{enumii}{13}
    	\item Group of Order 72. Eq. \eqref{eq: 5.13}.   
    	\item Group of Order 120. Eq. \eqref{eq: 5.15}.
    \end{enumerate}
\end{enumerate}

\section{p.22}

We now consider the case in which the normal curve of a structure of genus $p=4$ sits on a cone of the 2nd degree $K_2$. This special case has attracted the attention of Messrs. NOETHER and SCHOTTKY in the literature. As a strange property of the curve, it has been recognised that the same curve (in the typical form of a $C_9$ with 8 triple points) can be a coincidence curve with an unambiguous transformation of the plane in itself. One also has in itself a more general class of unambiguous transformations of the plane, that \guillemotright with 8 points \guillemotright, which is characterised by an unambiguous correspondence among the points \marginpar{start p.23} of such a $C_9$. For this reason Mr. S. KANTOR, in his investigations of the isolated types of unambiguous periodic transformations of the plane, looked for the cyclic transformations of the above structures with $p=4$ in themselves;\footnote{Acta Mathematica, Vol. 19, p. 164.} but not the complete groups, because that was not necessary for the problem to be solved. It is to this more general task that we shall concern ourselves here. We want to use a method that differs somewhat from KANTOR's. \\
If the curve is projected onto a plane from a point of the cone, the image curve is a $C_6$ with 2 infinitely close triple points, which come from the branches located on the generator of the cone through the starting point.\footnote{Cf. CLEBSCH-LINDEMANN, \textit{Vorlesungen} II, p. 431.} So the equation of the image curve can be brought to the following form: 
\[
z^3 y^3 + z^2 y^2 f_2(x,y) + zy f_4(x,y)+ f_6(x,y) = 0 \, , 
\]
where $y=x=0$ forms a triple point with $y=0$ as the common tangent of the 3 branches. The generators of the cone are mapped by the pencil of rays $x=\lambda y$, and the points of the curve occupied on one generator are represented by the equation: 
\[
z^3 + z^2 f_2(\lambda, 1) + z f_4(\lambda, 1) + f_6(\lambda, 1) = 0 = \varphi_3(z) \, , 
\] 
where $z=\infty$ is the vertex of the cone or the common point of all referred to as generators. The polars of the pole $z=\infty$ for $\varphi_3(z)$ are of course invariant structures. Here we consider the linear polar $3z + f_2(\lambda, 1)=0$ whose location satisfies the conic $3zy+f_2(x,y)=0$. But this conic section is the image of a cross-section of the cone. From now on we choose the centre of projection on this cross-section; the conic section then turns into a straight line, which we can denote by $z=0$, and for the general form of the equation of the image curve we get:
\begin{equation}\label{eq: 8.I}
z^3 y^3 + zy f_4(x,y) + f_6(x,y) = 0 \, . \tag{I}
\end{equation}
The collineation group of the normal curve occupied on the cone must be of such a nature that both a point\marginpar{start p.24}, the vertex of the cone, remains fixed as well as a plane, the location of the linear polars mentioned. The question arises whether, with a collineation belonging here, each generator can remain fixed; the 3 points of the curve must be cyclically permuted on a generator, so that the only possibility is a $G_3$. So we can state the theorem: \\
\textit{The collineation group of a normal curve of genus $p=4$ located on a cone $K_2$ is either a linear group or a combination of a linear group with a $G_3$.} \\
If we now start from the projected curve \eqref{eq: 8.I}, we have to consider a group of quadratic transformations instead of a collineation group. However, the distinguished straight line $z=0$ must always merge into itself, and 2 fundamental points must also lie in the infinitely neighboring triple points. The equation of a related transformation can be brought into the form:
\begin{equation}\label{eq: 8.a}
x^\prime : y^\prime : z^\prime = (a_1 x + b_1 y) (a_2 x + b_2 y) : (a_2 x + b_2 y)^2 : zy \, . \tag{a}
\end{equation}
However, it can be seen here that the substitution: 
\[
x^\prime : y^\prime = (a_1 x + b_1 y) : (a_2 x + b_2 y)
\] 
must convert the functions $f_6$ and $f_4$ into itself. One has therefore first to look for the common transformation group of these two functions. \\
The generators of the cone, which are fixed for a single spatial collineation, can be chosen for $x=0$ and $y=0$. The 4 fixed points in the collineation are the vertex of the cone and 3 points in the aforementioned fixed plane, 2 of which lie on those fixed generators, and the third is the pole of their connecting line with respect to the cross-section of the cone. But there are cases where all the points of a line joining two of these points, or of a plane through three of them, remain fixed. The planar transformation \eqref{eq: 8.a} turns into a collineation here:
\begin{equation}\label{eq: 8.aprime}
x^\prime : y^\prime : z^\prime = \alpha x : \beta y : \gamma z \, . \tag{$a^\prime$}
\end{equation}
This collineation transforms $f_6(x^\prime, y^\prime)$ into $k \cdot f_6(x, y)$ and $f_4(x^\prime, y^\prime)$ into $l \cdot f_4(x,y)$. The individual terms of \eqref{eq: 8.I}\marginpar{start p.25} are to be multiplied by the same number $k$. From this we get the conditions:
\begin{equation}\label{eq: 8.b}
\beta^3 \gamma^3 = \beta \gamma l = k \, ; \quad \gamma = \frac{k}{\beta l} \, ; \quad k^2 = l^3 \, . \tag{b} 
\end{equation}
Projecting the curve from the fixed point on the generator $x=0$, the equation of the image curve is: 
\[
z^3 x^3 + zx f_4(x, y) + f_6(x, y) = 0 \, , 
\] 
and the collineation is transformed into the following plane
\begin{equation}\label{eq: 8.a1prime}
x^\prime : y^\prime : z = \alpha x : \beta y : \gamma_1 z \, ; \quad \gamma_1 = \frac{k}{\alpha l} = \frac{\beta}{\alpha} \gamma \, . \tag{$a_1^\prime$}
\end{equation}
If one sets $\alpha=\beta$ in \eqref{eq: 8.aprime}, each generator is transformed into itself; one gets $k=\alpha^6$, $l=\alpha^4$ and $\gamma=\alpha$. The collineation is thus an identity. However, other conditions arise when $f_4$ vanishes identically, which, however, is initially excluded. \\
If $\beta = \gamma$, each point of the generator $x=0$ is transformed into itself. The same applies to the generator $y=0$ if $\alpha = \gamma_1$. If these relations are valid at the same time, one gets $\alpha^2 = \beta^2$ or $\alpha=-\beta$. The spatial collineation is then a harmonic homology whose plane of perspective goes through the solid generators. \\
If one takes $\alpha = \gamma$, then $\beta = \gamma_1$. The straight line with nothing but fixed points in the spatial collineation is here the polar line of the plane through the fixed generators. \\
By means of the above indications, the properties of the spatial collineations can be deduced by considerations in the plane.
\section{p.25}
First, we assume that $f_4$ does not vanish identically. The collineation group of the normal curve itself is obtained from the joint group of $f_4$ and $f_6$, which satisfies condition \eqref{eq: 8.b} $k^2=l^2$. The curve is touched by the 12 generators, $27 f_6^2 + 4 f_4^3=0$, of the cone. So that the curve does not have a double point, no 2 of those 12 straight lines may coincide; unless $f_6$ and $f_4$ share common factors that include stationary tangents; but $f_6$ must not contain a double or multiple factor, which also goes into $f_4$. Here is a brief list of the different types. \marginpar{start p.26}
\begin{enumerate}[label=\arabic*)]
    \item \label{item: 9.1}Perspective collineation of period 2: $x^\prime : y^\prime : z^\prime = -x : y : z$.
    \begin{equation}\label{eq: 9.1}
        z^3 y^3 + zy (ax^4 + b x^2 y^2 + c y^4) + dx^6 + e x^4 y^2 + f x^2 y^4 + g y^6 = 0 \, . \tag{1}
    \end{equation}
    \item \label{item: 9.2}Collineation \guillemotright with axes\guillemotright{} of period 2 : $x^\prime : y^\prime : z^\prime = x : -y : z$.
    \begin{equation}\label{eq: 9.2}
        z^3 y^2 + z(ax^4 + b x^2 y^2 + c y^4) + x(dx^4 + e x^2 y^2 + f y^4) = 0 \, . \tag{2}
    \end{equation}
    Here the order of the image curve is reduced because the centre of projection is a point on the curve. The axes are the straight lines connecting two fixed points of the curve in the invariant plane and their polar line.
    \item \label{item: 9.3}Klein four-group, consisting of 2 operations of type \ref{item: 9.1} and one of type \ref{item: 9.2}. One gets the equation: 
    \begin{equation}\label{eq: 9.3}
        z^3 y^2 + z[a(x^4 + y^4) + b x^2 y^2] + x[c(x^4 + y^4) + d x^2 y^2] = 0 \, . \tag{3}
    \end{equation}
    Here $f_6$ and $f_4$ merge into themselves when $x$ and $y$ are swapped.
    \item \label{item: 9.4}Klein four-group, which consists of nothing but operations of the type \ref{item: 9.2}. Since 2 points of $f_6=0$ remain fixed for each substitution, $f_6$ must be the JAKOBIan covariant of $f_4$.
    \begin{equation}\label{eq: 9.4}
        z^3 y^2 + z[a(x^4 + y^4) + bx^2 y^2] + x(x^4 - y^4) = 0 \, . \tag{4}
    \end{equation}
    \item \label{item: 9.5}Dihedral $G_8$, which contains as subgroups a Klein four-group of each of the two previous species.
    \begin{equation}\label{eq: 9.5}
        z^3 y^2 + a z x^2 y^2 + x(x^4 + y^4)= 0 \, . \tag{5}
    \end{equation}
    The distinguished cyclic subgroup derives from the repetition of the operation, $x^\prime y^\prime : z^\prime = x : iy : -z$.
\end{enumerate}
The 5 species thus obtained fully correspond in sequence to the 5 given in the 3rd section. The following species, however, is a genuinely new one
\begin{enumerate}[label=\arabic*)]
    \setcounter{enumi}{5}
    \item \label{item: 9.6}cyclic $G_4$ containing a $G_2$ of type \ref{item: 9.1} as a subgroup.
    \begin{equation}\label{eq: 9.6}
        z^3 y^2 + z(x^4 + a y^4) + y(bx^4 + c y^4) = 0 \, . \tag{6}
    \end{equation}
    $x^\prime : y^\prime : z^\prime = ix : y : z$.\marginpar{start p.27} \\
    The generator $x=0$ remains fixed at each point. In addition, one point of the curve moves into the projection centre on $y=0$, which causes the order of the image curve to decrease. The case follows from \eqref{eq: 9.1} if there $b=d=f=0$.
    \item \label{item: 9.7}Dihedral $G_6$. 
    \begin{equation}\label{eq: 9.7}
        z^3 y^3 + a z y^3 x^2 + x^6 + b x^3 y^3 + y^6 = 0 \, . \tag{7} 
    \end{equation}
    The generating operation of the cyclic $G_3$ is the following: $x^\prime : y^\prime : z^\prime = x : jy : z$. With this $G_3$ every point of the polar line remains fixed. The 3 $G_2$ are of type \ref{item: 9.1} and are obtained by swapping $x$ and $y$. If one has $b=0$ in \ref{item: 9.7}, the transformation group is a 
    \item \label{item: 9.8}Dihedral $G_{12}$. 
    \begin{equation}\label{eq: 9.8}
        z^3 y^3 + a z y^3 x^2 + x^6 + y^6 = 0 \, . \tag{8}
    \end{equation}
    The generating operation of the cyclic subgroup $G_6$ is: $x^\prime : y^\prime : z^\prime = -x : jy : z$. 
\end{enumerate}
Cases \ref{item: 9.7} and \ref{item: 9.8} correspond exactly to cases \eqref{eq: 4.6} and \eqref{eq: 4.7} of section 4.
\begin{enumerate}[label=\arabic*)]
    \setcounter{enumi}{8}
    \item \label{item: 9.9}Cyclic $G_3$ where the points of a generator remain fixed.
    \begin{equation}\label{eq: 9.9}
        z^3 y^3 + z y^2 (ax^3 + by^3) + x^6 + c x^3 y^3 + d y^6 = 0 \, . \tag{9}
    \end{equation}
    $x^\prime : y^\prime : z^\prime = jx : y : z$. 
\end{enumerate}
The points of $x=0$ remain fixed here. If $f_6$ is also the Jacobian covariant of $f_4$, the transformation group consists of 4 such $G_3$ and a Klein four-group of the kind \ref{item: 9.4}, which together form a
\begin{enumerate}[label=\arabic*)]
    \setcounter{enumi}{9}
    \item \label{item: 9.10}tetrahedral group. 
    \begin{equation}\label{eq: 9.10}
        z^3 y^3 + a z y^2(x^3 + y^3) + x^6 + 20x^3 y^3-8y^6 = 0 \, . \tag{10}
    \end{equation}
\end{enumerate}
One also gets the equation under the form \ref{item: 9.4} if only $b=2ai\sqrt{3}$, by which condition the equianharmonic invariant of $f_4$ vanishes. A special case of both (1) and (9) is the following.
\begin{enumerate}[label=\arabic*)]
    \setcounter{enumi}{10}
    \item \label{item: 9.11}Cyclic $G_6$.
    \begin{equation}\label{eq: 9.11}
        z^3 y^3 + a z y^5 + x^6 + b y^6 = 0 \, . \tag{11}
    \end{equation}
    $x^\prime : y^\prime : z^\prime = -jx : y: z$. \\
    \marginpar{start p.28} Also with this $G_6$ all points of the generator $x=0$ remain fixed. If one has $b=0$ here, the group is a 
	\item cyclic $G_{12}$.
	\begin{equation}\label{eq: 9.12}
	z^3 y^3 + zy^5 + x^6 = 0 \, . \tag{12}
	\end{equation}
	$x^\prime : y^\prime : z^\prime: = -jx : iy : -iz$. 
	\item Cyclic $G_5$.
	\begin{equation}\label{eq: 9.13}
	z^3 y^2 + azy^4 + x^5 + by^5 = 0 \, . \tag{13}
	\end{equation}
	$x^\prime : y^\prime : z^\prime = \varepsilon x : y: z$. \\
	The points of $x=0$ remain fixed. If $b=0$, the group is a 
	\item cyclic $G_{10}$. $x^\prime : y^\prime : z^\prime = \varepsilon x : -y : z$. 
	\begin{equation}\label{eq: 9.14}
	z^3 y^2 + zy^4 + x^5 = 0 \, . \tag{14}
	\end{equation}
	The subgroup $G_2$ here is of the type \ref{item: 9.2}. 
\end{enumerate}
In the remaining cases, $f_4$ vanishes identically. The normal curve always merges into itself through a perspective $G_3$, which leaves each generator unchanged. In general, 
\begin{enumerate}[label=\arabic*)]
    \setcounter{enumi}{14}
    \item this $G_3$, $x^\prime: y^\prime : z^\prime = x:y:jz$, is also the complete group of the curve 
    \begin{equation}\label{eq: 9.15}
    z^3 y^3 = f_6(x, y) \, . \tag{15}
    \end{equation}
\end{enumerate}
But if $f_6$ has a linear transformation group in itself of order $\mu$, one obtains the group of the curve by combining this $G_\mu$ with $G_3$, from which the order $3\mu$ results. Within this $G_{3 \mu}$ the perspective $G_3$ is normal. \\
The possible groups of $f_6$ in itself have already been compiled in our earlier discussion.\footnote{Bihang till K. Sv. Vet.-Akad. Handl. Vol 21, Section. I, No. 1} We found there in the 4th section that the same are the following: $G_2$, Klein four-group, Dihedral $G_6$, $G_5$, Dihedral $G_{12}$ and octahedral group. We can therefore immediately write down the following normal equations:
\begin{align}
	z^3 y^3 &= x^6 + ax^4 y^2 + bx^2 y^4 + y^6 \, ; \label{eq: 9.16} \tag{16} \\
	z^3 y^2 &= x(x^4 + ax^2 y^2 + y^4) \, ; \label{eq: 9.17} \tag{17} \\ 
	z^3 y^3 &= x^6 + a x^3 y^3 + y^6 \, ; \label{eq: 9.18} \tag{18} \\
	z^3 y^2 &= x^5 + y^5 \, ; \label{eq: 9.19} \tag{19} \\ 
	z^3 y^3 &= x^6 + y^6 \, ; \label{eq: 9.20} \tag{20} \\ 
	z^3 y^2 &= x(x^4 + y^4) \, . \label{eq: 9.21} \tag{21} 
\end{align}\marginpar{start p.29}
The group of \eqref{eq: 9.16} is a cyclic $G_6 : x^\prime: y^\prime : z^\prime = - x : y: jz$, and that of \eqref{eq: 9.19} is a cyclic $G_{15} : x ^\prime : y^\prime : z^\prime = \varepsilon x : y : jz$. The group of \eqref{eq: 9.17} is of order 12 and consists of 3 cyclic $G_6$, all of which have the perspective $G_3$ as a subgroup. Of these $G_6$, 2 are of the same nature as that of Curve \eqref{eq: 9.16}; but the third, $x^\prime : y^\prime : z^\prime = x : -y : jz$, has a $G_2$ \guillemotright with axes\guillemotright{} as a subgroup. We do not dwell on the groups of curves \eqref{eq: 9.18} and \eqref{eq: 9.20} of the respective orders 18 and 36. The group of order 72 of curve \eqref{eq: 9.21} has 3 cyclic $G_{12}$, which are assigned to the cyclic $G_4$ of an octahedron group. The one $G_{12}$ is transformed into the plane collineation, $x^\prime : y^\prime : z^\prime = x : iy : -jz$. \\
Regarding the cyclic groups, our results differ in a few points from those given by Mr. KANTOR.\footnote{Ibid, p. 166.} Mr. KANTOR does not mention the $G_2$ \guillemotright with axes\guillemotright{} (our case 2) and the cyclic $G_6$ belonging to the curve \eqref{eq: 9.17}, which has a perspective $G_3$ and a $G_2$ \guillemotright with axes\guillemotright{} as subgroups. We have found 2 types of cyclic $G_{12}$ here: one belongs to the curve \eqref{eq: 9.12} and contains as subgroups a perspective $G_2$ and a $G_3$ with a generator fixed at each point, the other to the curve \eqref{eq: 9.21} and contains a perspective $G_3$ and a $G_2$ \guillemotright with axes\guillemotright. Mr. KANTOR considers these two cases (designated 12 and 16 in his account) to be equivalent, which is not necessarily true. On the other hand, KANTOR's forms 6 and 15, 9 and 10 are equivalent to each other (and to our cases 6 and 9), but it should be noted that KANTOR's forms 6, 15 and 10 are not given with the full number of terms. 
\section{p.29}
We now turn to the treatment of structures of genus $p=5$. The plane curve of the lowest order, into which a general structure of this genus can be unambiguously transformed, is known to be a $C_6$ with 5 double points; but if the structure has a special linear series $g_3^1$, it can be transformed into a $C_5$ with a double point.\footnote{See CLEBSCH-LINDEMANN, \textit{Vorlesungen \"uber Geometrie}, I, p. 709.}\marginpar{start p.30} In the latter case, the pencil of rays intersects the $g_3^1$ through the double point, and the adjoint $\varphi$-curves are conics through the same point. If one now goes to the $\varphi$-normal-curve in 4-dimensional space, one finds that it must lie on a third-order ruled surface, which is the image of our original plane, with the rays from the double point being transformed into generators and the double point itself in be transformed into a simple straight line of the ruled surface. Each generator is thus a trisecant and the directrix a bisecant of the normal curve. From this it can be seen that the ruled surface or the plane of the $C_5$ plays an distinguished role. Consequently, the unambiguous transformations of the $C_5$ must be birational in themselves, because no quadratic or higher collineations are possible. The latter is concluded from the fact that the system of adjoint conics must also be transformed into itself in such a birational transformation; but within this system there is an distinguished degenerating subset, and the associated conics each consist of a straight line through the double point and an arbitrary straight line; the system of straight lines must therefore be converted into itself, which only occurs with a collineation. The possible collision groups of $C_5$ must be of a very trivial nature. One descends from the total group to the normal subgroup, in which the two branches in the double point remain uninterchanged, and from this to a (cyclic) normal subgroup, in which each straight line merges into itself through the double point. However, these normal groups may consist of the total group or just identity. \\
According to a theorem by Mr. WEBER, an algebraic structure with any $p$ is characterised by $\frac{1}{2}(p-2)(p-3)$ linearly independent quadratic relations between the $\varphi$-functions.\footnote{Math. Ann., Vol. 13.} However, these relations are not sufficient to define the structure if a $g_3^1$ exists.\footnote{Based on a theorem of L. KRAUS, Math. Ann., Vol. 16. It should be noted here that the relations $\Phi^{(2)}=0$ are not sufficient to define the structure even in a case without $g_3^1$, namely for the plane $C_5$ without double points.} So\marginpar{start p.31} for $p=5$ we have three quadratic relations, $\Phi_1^{(2)}=0$, $\Phi_2^{(2)}=0$, $\Phi_3^{(2)}=0 $, the normal curve $C_8$ can generally be considered as its intersection curve;\footnote{Cf. a treatise by Mr. NOETHER, Math. Ann., Vol. 26.} but the $\Phi^{(2)}=0$ in the case of the $g_3^1$ have a common surface, namely the previously mentioned ruled surface of the third order.
\section{p.31}
From now on we stick to the general case. The normal curve is determined in the 4-dimensional space by the three quadratic relations: 
\begin{equation}\label{eq: 11.1}
F_1(x_1, x_2, x_3, x_4, x_5) = 0 \, , \quad F_2=0 \, , \quad F_3=0 \, . \tag{1}
\end{equation}
We take any linear combination of $F_i$: 
\begin{equation}\label{eq: 11.2}
\lambda_1 F_1 + \lambda_2 F_2 + \lambda_3 F_3 = 0 \tag{2}
\end{equation} 
and look for its discriminant: 
\begin{equation}\label{eq: 11.3}
\Delta_5(\lambda_1, \lambda_2, \lambda_3) = 0 \, , \tag{3}
\end{equation} 
which can be written in the known way in the form of a determinant with 5 rows and 5 columns. \textit{The nature of the covariant curve $\Delta_5=0$ in the $\lambda$ plane is of fundamental importance for our normal curve.} Equation \eqref{eq: 11.3} determines those combinations of $\lambda_i$ for which the left term of \eqref{eq: 11.2} is a function of only 4 variables, say $x_1^\prime$, $x_2^\prime$, $x_3^\prime$, $x_4^\prime$. The location of the vertices $x_1^\prime = x_2^\prime = x_3^\prime = x_4^\prime=0$ clearly corresponds to the curve $\Delta_5$.
However, this location is determined by the determinant group 
\begin{equation}\label{eq: 11.4}
\begin{vmatrix}
F_1^1 & F_1^2 & F_1^3 & F_1^4 & F_1^5 \\ F_2^1 & F_2^2 & F_2^3 & F_2^4 & F_2^5 \\ F_3^1 & F_3^2 & F_3^3 & F_3^4 & F_3^5 
\end{vmatrix} = 0 \, , \tag{4}
\end{equation}
where $F_i^k = \frac{\partial F_i}{\partial x_k}$. The system \eqref{eq: 11.4} represents a curve of order 10 in the 4-dimensional space,\footnote{See SALMON-FIEDLER, \textit{Algebra der linearen Transformationen}, 2. Edition (1877), p. 368.} which we will denote by $D_{10}$. The normal curve cannot have a double point because then its genus is reduced.\marginpar{start p.32} However, as is easy to find, both equations \eqref{eq: 11.1} and \eqref{eq: 11.4} should hold for a double point. The necessary condition follows that the curve $D_{10}$ must not separate the normal curve. Thus, no double points of $D_{10}$ and $\Delta_5$ can come from such intersections. But the curve $\Delta_5$ can have double points, namely in those points of the $\lambda$-plane, for which the left term of \eqref{eq: 11.2} as a function of only 3 variables, say $x_1^\prime$, $x_2^\prime$, $x_3^\prime$, can be written; such a double point corresponds to a straight line on $D_{10}$, $x_1^\prime=x_2^\prime=x_3^\prime=0$. The relations $F_i=0$ between 3 variables correspond to just as many systems of \textit{infinitely many ABELian root forms}. L. KRAUS now wanted to show that there are at most 3 such systems;\footnote{Ibid.} but his proof is incorrect, and we shall find that as many as 10 such systems can occur. In this case, $\Delta_5$ has 10 double points, which is only possible if it consists of 5 lines, and $D_{10}$ consists of 10 lines, for which each 4 intersect at 5 points, these lines and multiple points of $D_{10}$ correspond to double points and straight lines of $\Delta_5$. But we shall come back to this case later. \\
According to BRILL-NOETHER's law of reciprocity, the special linear series $g_4^1$ of the structures of genus $p=5$ are assigned to each other in pairs. The curve $D_{10}$ has a simple connection with the system of $g_4^1$. If one projects the normal curve from a point of $D_{10}$, one obtains a curve $C_8$ as its image in 3-dimensional space, which extends to a surface of the 2nd degree $F(x_1^\prime, x_2^ \prime, x_3^\prime, x_4^\prime)=0$, where $F=0$ gives the corresponding relation with only 4 variables. The generating systems of this surface cut out 2 reciprocal $g_4^1$ to each other. But if that relation were between only 3 antecedents, the projected curve lies on a second-degree cone, so that the two series $g_4^1$ coincide, and the only $g_4^1$ is its own reciprocal; the infinitely many associated ABELian root forms are assigned to the tangent planes of the cone. Thus the system of $g_4^1$ is 1-2-unambiguously related to the curve $D_{10}$, as we later want to illustrate geometrically in a different way.\marginpar{start p.33} With a collineation of the normal curve in itself, the system of relations $\lambda_1 F_1 + \lambda_2 F_2 + \lambda_3 F_3=0$ must of course also be converted into itself. The $\lambda_i$ thus experience linear substitutions, in which the covariant curve $\Delta_5$ is necessarily transformed into itself. This is not contradicted by the fact that a collineation of the $\Delta_5$ is not always the result of a collineation of the normal curve. On the other hand, the normal curve can undergo unambiguous transformations in itself in which each point of $\Delta_5$ remains fixed, and therefore each relation \eqref{eq: 11.2} is retained unchanged. With such a collineation, the curve $D_{10}$ must also merge point by point, if one disregards straight lines that correspond to double points of $\Delta_5$; the remaining curve $D_{10}$ can consequently not be an actual 4-dimensional curve, but can either be composed of a point and a 3-dimensional curve or of a straight line and a flat curve, depending on a point in the collineation and the points of a space or the points of a straight line and a plane remain fixed. As a necessary condition for this one immediately obtains that each relation \eqref{eq: 11.1} $F_i=0$ can be written either in the form 
\begin{equation}\label{eq: 11.a}
a_i x_1^2 + f_i(x_2, x_3, x_4, x_5) = 0 \quad (i=1,2,3) \tag{a}
\end{equation} 
or in the form 
\begin{equation}\label{eq: 11.b}
\varphi_i(x_1, x_2) + \psi_i(x_3, x_4, x_5) = 0 \, . \tag{b}
\end{equation} 
In the first case one has the collineation: 
\begin{equation}\label{eq: 11.alpha}
x_1^\prime = -x_1 \, , \quad x_i^\prime = x_i \quad (i=2,3,4,5) \tag{$\alpha$}
\end{equation} 
with 8 fixed points of the curve in space $x_1=0$; in the latter case the collineation: 
\begin{equation}\label{eq: 11.beta}
x_1^\prime = -x_1 \, , \quad x_2^\prime = -x_2 \, , \quad x_i^\prime = x_i \quad (i=3,4,5) \tag{$\beta$}
\end{equation} 
without fixed points of the curve. The straight lines connecting corresponding points in cases \eqref{eq: 11.a} and \eqref{eq: 11.b} thus form (according to equation \eqref{eq: A} introduction) a ruled surface of genus $p^\prime=1$ or $p^\prime=3$. However, symmetric transformations with the assigned genus $p^\prime=2$ can also occur, but without transforming every point of the curve $\Delta_5$ into itself. In case \eqref{eq: 11.a} the discriminant of equation \eqref{eq: 11.2} breaks up into a linear and a biquadratic \marginpar{start p.34}factor and in case \eqref{eq: 11.b} into a quadratic and a cubic factor, so that the curve $\Delta_5$ consists of a straight line and a $C_4$, respectively, a conic section and a $C_3$. But these curves can also break up, so that the following 7 typical cases are obtained.
\begin{enumerate}[label=\arabic*)]
	\item $\Delta_5$ is an actual $C_3$.
	\item $\Delta_5$ breaks up into a straight line and a $C_4$. The intersection points of the straight lines with $C_4$ correspond to 4 systems of infinitely many ABELian root forms. But here, as in the previous case, other such systems can appear, corresponding to the possible double points of $C_4$ or $C_5$.
	\item $\Delta_5$ breaks up into a conic section and a $C_3$. The 6 points of intersection of these curves correspond to 6 systems of ABELian root forms.
	\item $\Delta_5$ breaks up into 2 straight lines and a $C_3$. Any relation $F_i=0$ has the form
	\begin{equation}\label{eq: 11.c}
	a_i x_1^2 + b_i x_2^2 + \psi_i(x_3, x_4, x_5) = 0 \, . \tag{c}
	\end{equation}
	Here one has 2 collineations of type \eqref{eq: 11.alpha} and one of type \eqref{eq: 11.beta}. Seven systems of ABELian root forms always occur.
	\item $\Delta_5$ breaks up into a straight line and 2 conics.
	\begin{equation}\label{eq: 11.d}
	    F_i = a_i x_1^2 + \varphi_i(x_2, x_3) + \psi_i(x_4, x_5) = 0 \, . \tag{d}
	\end{equation}
	A collineation of type \eqref{eq: 11.alpha} and 2 of type \eqref{eq: 11.beta}. Eight systems of ABELian root forms.
	\item $\Delta_5$ breaks up into 3 lines and a conic section.
	\begin{equation}\label{eq: 11.e}
	    F_i = a_i x_1^2 + b_i x_2^2 + c_i x_3^2 + \psi_i(x_4, x_5) = 0 \, . \tag{e}
	\end{equation}
	There are 3 collineations of type \eqref{eq: 11.alpha} and 4 of type \eqref{eq: 11.beta}. Nine systems of ABELian root forms.
	\item $\Delta_5$ breaks up into 5 straight lines. 
	\begin{equation}\label{eq: 11.f}
	F_i = a_i x_1^2 + b_i x_2^2 + c_i x_3^3 + d_i x_4^4 + e_i x_5^2 = 0 \, . \tag{f}
	\end{equation}
	In this case one has 5 collineations of type \eqref{eq: 11.alpha} and 10 of type \eqref{eq: 11.beta}. From the 3 independent relations $F_i=0$ one can eliminate any pair $x_k^2$ and $x_i^2$ and thus obtain 10 relations between 3 variables, which correspond to just as many systems of Abel root forms.
\end{enumerate}
\marginpar{start p.35} The question can now be raised as to how one can arrive at the structure in 4-dimensional space given $\Delta_5$. First it should be emphasised that both a general $C_5$ and a structure of the genus $p=5$ have twelve moduli that are indestructible by unambiguous transformation, so that a finite number of solutions can be expected. It turned out
\[
\Delta_5(\lambda) = \sum \pm u_{11} u_{22} u_{33} u_{44} u_{55}
\] 
in determinant form, where the $u_{ik}=u_{kl}$ are linear functions of the $\lambda$. The unambiguous transformation connecting the curves $\Delta_5$ and $D_{10}$ has the form 
\[
x_1 : x_2 : x_3 : x_4 : x_5 = U_{i1} : U_{i2} : U_{i3} : U_{i4} : U_{i5}
\] 
or 
\[
x_i x_k = \rho U_{ik} \, , 
\] 
where $U_{ik}$ is a subdeterminant from $\Delta_5$ to $u_{ik} $. From these equations one can produce the $\lambda_1^4$, $\lambda_1^3 \lambda_2, \, \dots$ as a linear function of the $x_1^2$, $x_1 x_2, \, \dots$. According to this, the intersections of $\Delta_5$ with all fourth-order curves correspond to the intersections of $D_{10}$ with the manifolds represented by the quadratic relations; \textit{in particular, the intersections of $D_{10}$ with the doubly counted 3-dimensional spaces correspond to a quadruple infinite system of fourth-order tangent curves at $\Delta_5$}. This system results from changing the determinant $\Delta_5$. So one would have to look for that class of touch [contact. ID: tangent?, LDH: Maybe, discuss further with AB] systems by which $\Delta_5$ can be written in the form of a symmetric determinant.\footnote{Cf. HESSE's investigations on the contact conic sections of a $C_3$ and on the systems of the \textit{first} kind of $C_3$ in contact with a $C_4$, CRELLE's Journal, Vol. 36,49.} The formations in the 4-dimensional space originating from different contact systems are generally not collinearly related, because for this it would be necessary for the systems to also be related in this way. It should also be noted that a collineation of the $\Delta_5$ must also convert the underlying contact system into itself, so that a corresponding collineation exists for the normal curve with $p=5$. The contact systems in question must be adjoint \marginpar{start p.35}and pass through the double points of $\Delta_5$, which is evident from the fact that the double points of the $\Delta_5$ correspond to straight lines of $D_{10}$.
\section{p.36}
A $C_7$ in the usual 3-dimensional space and a flat $C_6$ can always, if their genus $p=5$, come out as projections of the normal curve, because the $\infty^2$ planes resp. $\infty^2$ straight lines cut out special linear series in the two cases mentioned. Now the question arises how, given a plane $C_6$ with 5 double points, one can determine the two points of the curve which are to be regarded as successive projection centres. These points, together with an arbitrary pair united in a double point, should form a group $G_4$, which belongs to a special linear series $g_4^1$. However, this condition is only satisfied by those points which are cut out by a conic section drawn through the 5 double points. \\
The $C_7$ projected from the various points of the normal curve each reveal the peculiarities of the system of special linear series $g_4^1$ through their \textit{trisecant ruled surface}. The 3 points which, together with the projection centre, form a group $G_4$ of such a special linear series must be projected in a straight line. \textit{Hence the generators of the ruled surface and the series $g_4^1$ correspond unambiguously.} The point groups of the reciprocal $g_4^1$ are cut from the planes through the trisecant associated to the given $g_4^1$, and in particular 2 generators corresponding to reciprocal $g_4^1$ intersect each other, their common plane going through the image point of the projection centre. The locus of the intersections of two such reciprocal generators is the projected image of the curve $D_{10}$. If a $g_4^1$ and consequently the corresponding trisecant is reciprocal to itself, the trisecant must be double torsal such that the torsal plane contains the tangents of the 3 points lying on the torsal of the $C_7$.
The $C_7$ is the fivefold curve and the $D_{10}$ double curve of the trisecant end ruled surface. A generator thus intersects 12 others in the 3 points of $C_7$ and one in the point of $D_{10}$, and the degree of the ruled surface is consequently 15. For the genus $P$ of a cross-section one obtains according to PLUCKER's formulas: 
\[
P=11-\delta \, , 
\]
\marginpar{start p.37}where $\delta$ means the number of the above-mentioned double torsals with the same plane. Depending on the nature of $\Delta_5$ and $D_{10}$ one has (corresponding to the typical cases of the previous section) the following 7 cases.
\begin{enumerate}[label=\arabic*)]
	\item $\Delta_5$ is an actual $C_5$, which can have a maximum of 6 double points. The genus of $R_{15}$ is usually 11, but can be lowered to 5 by $\delta$.
	\item $\Delta_5 = C_4 + C_1$. $R_{15} = R_{12} + K_3$. The genus of $R_{12}$ alternates between 7 and 4, with the reduction being effected by double points of $C_4$. $K_3$ is of genus 1, and the cusp contains that point of $C_7$ which corresponds to the centre of projection in the transformation \eqref{eq: 11.alpha} of the previous section. The generators of the cone connect the points corresponding to this transformation. $K_3$ and $R_{12}$ intersect in $C_7$, counted quadruple, touching each other along the 4 double torsals, which correspond to the intersections of the constituent parts of $\Delta_5$. Also in the following cases, in which $R_{15}$ breaks up, the sectional curves only consist of $C_7$ and double torsals.
	\item $\Delta_5 = C_3 + C_2$. $R_{15} = R_9 + R_6$. $R_9$ has genus 4 or (through a double point of $C_3$) 3. $R_6$ is always of genus $P=2$.\footnote{Regarding this $R_6$ see my work, \textit{Klassifikation af regelytorna af 6. graden} p. 57 (Lund 1892. Diss.)} 
	\item $\Delta_5 = C_3 + 2 C_1$. $R_{15} = R_9 + 2 K_3$. 
	\item $\Delta_5 = 2 C_2 + C_1$. $R_{15} = 2 R_6 + K_3$.
	\item $\Delta_5 = C_2 + 3 C_1$. $R_{15} = R_6 + 3 K_3$.
	\item $\Delta_5 = 5 C_1$. $R_{15} = 5 K_3$. The $C_7$ lies here on $5K_3$, each pair of which touches each other along the connecting line of the cusps.
\end{enumerate}
\section{p.37}
In order to determine the complete collineation group of a given normal curve, one first has to consider the covariant curve $\Delta_5$. One classifies them in one of the 7 cases of the 11th section and then look for their collineation group. From this group one has to remove the subgroup to which collineations of the normal curve also belong. A given linear substitution in the $\lambda_i$ corresponds to the transposed substitution in the $F_i$ due to the equation $\lambda_1 F_1 + \lambda_2 F_2 + \lambda_3 F_3 = 0$. Consequently, one can search for the typical periodic collineations in the 4-dimensional space from the outset, through which a linear transformation is effected under \marginpar{start p.38}3 quadratic relations whose intersection curve \textit{has no double point}. If we transfer the result of such an investigation to the $\lambda$-plane, it follows that only the collineations from period 2 and the non-perspective $G_3$, $G_4$ and $G_5$ need to be considered. In 4-dimensional space there are also cyclic $G_6$ and $G_8$; however, the associated subgroup $G_2$ leaves each quadratic relation unchanged. \\
The finite planar and \textit{simple} collineation groups are (apart from the cyclic ones with a prime number period) the icosahedron group, the $G_{168}$ known from the theory of a special $C_4$ and a $G_{360}$.\footnote{This $G_{360}$ was discovered by Mr. VALENTINE (Kong. Dansk. Vid. Selsk. Skrift. (6) V, 1889), but otherwise seems to be unknown.} The first two are known not to have an invariant form of 5th order; the same must be true for $G_{360}$ because it contains icosahedral groups. If one then also takes the results of section 10 regarding the structures with a $g_3^1$ into account, we have the following theorem: \\
\textit{All groups occurring within the non-hyperelliptic structures of genus $p=5$ can be broken up into a series of purely cyclic groups.} \\
If $\Delta_5$ now has a group of order $\mu$, to which collineations of the normal curve also belong, then one has the following numbers for the order of the total group of the latter, depending on the typical case of $\Delta_5$: $\mu, \, 2\mu, \, 2\mu, \, 4\mu, \, 4\mu, \, 8\mu, \, 16\mu$. \\
Finally, those curves are listed here in their normal equations whose moduli are completely determined by the condition of admitting certain collineation groups.
\begin{enumerate}[label=\arabic*)]
	\item \textit{Group of Order 192}.
	\begin{align*}
		F_1 &= x_1^2 + x_4^2 + x_5^2 = 0 \, , \\
		F_2 &= x_2^2 + x_4^2 - x_5^2 = 0 \, , \\
		F_3 &= x_3^2 + x_4 x_5 = 0 \, . 
	\end{align*}
    One has (with the identity) 8 collineations by sign changes from the $x_i$, where each $F_i$ is unchanged. Then one has to notice that the form $x_4 x_3 \cdot (x_4^2 + x_5^2) \cdot (x_4^2 - x_5^2)$ is an octahedron form broken up into its 3 distinct quadratic factors\marginpar{start p.39}. One uses the group of this binary form with suitable permutations and substitutions of $x_1, \, x_2, \, x_3$ and then by combining it with $G_8$ the complete $G_{192}$ is obtained.\footnote{In the theory of elliptic modular functions, this $G_{192}$ belongs to the eighth order main congruence group. See KLEIN-FRICKE, \textit{Modulfunctionen} I, p. 338ff.} Eliminating $x_5$ gives the equations:
    \[
    x_1^2 + x_2^2 + 2x_4^2 = 0 \, , \quad x_4^2(x_2^2 + x_4^2) - x_3^4 = 0 \, . 
    \]
    The curve $\Delta_5$ consists of a conic section and 3 straight lines, which form a polar triangle in relation to the conic section. In the 3 following cases, $\Delta_5$ consists of 5 straight lines, which are also in a special position.
    \item \textit{Group of order 64.}
    \begin{align*}
    	x_1^2 + x_2^2 + x_3^2 + x_4^2 + x_5^2 &= 0 \, , \\
    	x_1^2 + i x_2^2 - x_3^2 - i x_4^2 &= 0 \, , \quad (i^2=-1.) \\
    	x_1^2 - x_2^2 + x_3^3 - x_4^2 &= 0 \, . 
    \end{align*}
    The $G_{64}$ is created by combining the sign changes with the cyclic permutation of $x_1, ,\ x_2, \, x_3, \, x_4.$
    \item \textit{Group of order 96.}
    \begin{align*}
    	x_1^2 + x_4^2 + x_5^2 &= 0 \, , \\
    	x_2^2 + j x_4^2 + j^2 x_5^2 &= 0 \, , \quad (j^3=1.) \\
    	x_3^2 + j^2 x_4^2 + j x_5^2 &= 0 \, . 
    \end{align*}
    The $G_{96}$ is obtained by combining the sign changes and swapping $x_1, \, x_2$ and $x_3$ on the one hand, and $x_4$ and $x_5$ on the other.
    \item \textit{Group of order 160.}
    \begin{align*}
    	x_1^2 + x_2^2 + x_3^2 + x_4^2 + x_5^2 &= 0 \, , \\
    	x_1^2 + \varepsilon x_2^2 + \varepsilon^2 x_3^2 + \varepsilon^3 x_4^2 + \varepsilon^4 x_5^2 &= 0 \, , \quad (\varepsilon^5=1.) \\
    	\varepsilon^4 x_1^2 + \varepsilon^3 x_2^2 + \varepsilon^2 x_3^2 + \varepsilon x_4^2 + x_5^2 &=0 \, . 
    \end{align*}
    \marginpar{start p.40}Here the allowed permutations of the $x_i$ form a Dihedral $G_{10}$. Namely, the array $(x_1 \, x_2 \ x_3 \, x_4 \, x_5)$ can both be cyclically shifted into itself and reversed into the array $(x_5 \, x_4 \, x_3 \, x_2 \, x_1)$. The $G_{160}$ is created by combining this $G_{10}$ with the sign changes.
\end{enumerate}

\section{p.40}
Finally, a few hints about the formations of genus $p=6$ shall be permitted. Apart from the hyperelliptic case, there are 3 essentially different types, which we want to treat here in series.
\begin{enumerate}[label=\arabic*)]
	\item The curve has a special linear series $g_5^2$, i.e. it can be transformed into a plane $C_5$. With an unambiguous transformation of $C_5$ in itself, the point groups $G_5$ of $g_5^2$ must form a closed invariant system. Such a transformation must consequently be a collineation, and according to a theorem in the previous section composed entirely of cyclic groups.
	\item The curve has a special linear series $g_3^1$. Such a curve is easily unambiguously converted into a $C_6$ with a triple and a double point. In the case of a group of unambiguous transformations of the curve in itself, the point groups $G_3$ of the special linear series must be permuted in a closed manner, which of course can only be effected by one of the finite linear groups; but not by an icosahedron group, because the number of $G_3$ containing 2 coincident points is 16, and no 16 points in this group form a closed system. In addition, the transformation group can still be composed of a cyclic group of period 3, which cyclically exchanges the points of each $G_3$.
	\item In the general case, the curve has 5 $g_4^1$ and can be unambiguously converted into a $C_6$ with 4 double points, whereby the point groups of the $g_4^1$ are cut out by the pencil of rays by a double point each and the conic pencil by all 4 $\delta$. With an unambiguous transformation of the curve in itself, these $g_4^1$ can either be permuted in themselves or in some way. The former is not possible if there are no 3 double points in a straight line. Each permutation of the $g_4^1$ can be effected by collineations and successive quadratic transformations. One also gets a case where all 120 permutations of the 5 $g_4^1$ transform the curve into itself:\marginpar{start p.41} 
	\[
	2 \sum (x^4 y z + y^3 z^3) - 2 \sum x^4 y^2 + \sum x^3 y^2 z - 6 x^2 y^2 z^2 = 0 \, . 
	\]
	The double points are here in the coordinate corners and in the point $x=y=z$.
\end{enumerate}
Finally, as a result of our investigations, we can cite the following general propositions, which Mr. DYCK has already given for the genera $p=0, \, 1, \, 2, \, 3$.\footnote{\textit{\"Uber regul\"are Riemann'sche Fl\"achen}, Math. Ann., Vol. 17, p. 474} \\
\textit{The only simple and finite groups of unambiguous transformations of a curve in itself within the genera $p=0, \dots, 6$ are the cyclic groups of prime order, the icosahedron group and a $G_{168}$ (at $p=3$). The composite groups within our genera, with 3 exceptions, can be broken up into a sequence of merely cyclic groups. The exceptions are: BRING's curve ($p=4$) and the just mentioned curve with $p=6$, whose groups are holohedrally isomorphic with the group of 120 permutations of 5 things, as well as a hyperelliptic curve\footnote{The equation of this curve is $y^2 = x(x^{10} + 11x^5-1)$. In the theory of elliptic modular functions, it supplies the modules belonging to a normal congruence subgroup of the 10th order. See KLEIN-FRICKE, \textit{Modulfunctionen} I, p. 651} of genus $p=5$, whose group is also of order 120; namely, within these groups, an icosahedron group is normal.}

\end{document}